\documentclass[aop,citesort,MSNbibl,seceqn,dvips]{arximspdf}

\doi{10.1214/08-AOP402}
\volume{37}
\issue{1}
\pubyear{2009}
\firstpage{250}
\lastpage{274}

\makeatletter

\newtheorem{proposition}{Proposition}[section]
\newtheorem{theorem}[proposition]{Theorem}
\newtheorem{lemma}[proposition]{Lemma}
\newproclaim{example}[proposition]{Example}

\makeatother

\begin{document}
\begin{frontmatter}

\title{Continuity properties and infinite divisibility
of~stationary distributions of some generalized
Ornstein--Uhlenbeck processes}
\runtitle{Properties of stationary law of generalized O--U processes}

\begin{aug}
\author[A]{\fnms{Alexander} \snm{Lindner}\corref{}\ead[label=e1]{a.lindner@tu-bs.de}} and
\author[B]{\fnms{Ken-iti} \snm{Sato}\ead[label=e2]{ken-iti.sato@nifty.ne.jp}}
\runauthor{A. Lindner and K. Sato}
\affiliation{Technische Universit\"at of Braunschweig and
Hachiman-yama, Nagoya}
\address[A]{Institute for Mathematical Stochastics\\
TU Braunschweig\\
Pockelsstra{\ss}e 14\\
D-38106 Braunschweig\\
Germany\\
\printead{e1}} 
\address[B]{Hachiman-yama 1101-5-103\\
Tenpaku-ku\\
Nagoya 468-0074\\
Japan\\
\printead{e2}}
\end{aug}

\received{\smonth{2} \syear{2007}}
\revised{\smonth{2} \syear{2008}}

%
\begin{abstract}
Properties of the law $\mu$ of
the integral
$\int_0^{\infty}c^{- N_{t-}} \, dY_t$ are studied, where
$c>1$ and $\{ (N_t,Y_t), t\geq0\}$ is a bivariate L\'evy
process such that $\{ N_t \}$ and $\{ Y_t \}$ are Poisson processes
with parameters $a$ and $b$, respectively.
This is the stationary
distribution of some generalized Ornstein--Uhlenbeck process.
The law $\mu$ is parametrized by $c$, $q$ and $r$, where
$p=1-q-r$, $q$, and $r$ are the normalized L\'evy measure of
$\{ (N_t,Y_t)\}$ at the points $(1,0)$, $(0,1)$ and $(1,1)$,
respectively.
It is shown that, under the condition that $p>0$ and $q>0$,
$\mu_{c,q,r}$ is infinitely divisible if and only
if $r\leq pq$. The infinite divisibility of the symmetrization
of $\mu$ is also characterized. The law $\mu$ is
either continuous-singular or
absolutely continuous, unless $r=1$.
It is shown that if $c$ is in the set of
Pisot--Vijayaraghavan numbers, which
includes all integers bigger than~$1$, then
$\mu$ is continuous-singular under the condition $q>0$.
On the other hand, for Lebesgue almost every
$c>1$, there are positive constants $C_1$~and~$C_2$ such that $\mu$
is absolutely continuous whenever $q\geq C_1p \geq C_2 r$.
For any $c>1$ there is a positive constant $C_3$
such that $\mu$
is continuous-singular whenever $q>0$ and $\max\{q,r\}\leq C_3 p$.
Here, if $\{ N_t \}$ and $\{ Y_t \}$ are independent,
then $r=0$ and $q=b/(a+b)$.
\end{abstract}

%
\pdfkeywords{Decomposable distribution,
generalized Ornstein--Uhlenbeck process,
infinite divisibility,
Levy process,
Peres--Solomyak (P.S.) number,
Pisot--Vijayaraghavan (P.V.) number,
symmetrization of distribution}

\end{frontmatter}

\section{Introduction}\label{sec1}
A generalized Ornstein--Uhlenbeck process $\{ V_t, t\geq0\}$
with initial condition $V_0$ is defined as
\[
V_t = e^{-\xi_t} \biggl( V_0 + \int_0^t e^{\xi_{s-}} \, d\eta_s
\biggr),
\]
where $\{ (\xi_t, \eta_t), t \geq0\}$ is a bivariate L\'evy
process, independent of $V_0$. See Carmona, Petit and Yor \cite{CPY97,CPY01}
for basic properties.
Such processes arise in a variety of
situations such as risk theory (e.g., Paulsen \cite{Pa98}),
option pricing (e.g., Yor \cite{Yo}) or financial time series
(e.g., Kl\"uppelberg, Lindner and Maller \cite{KLM06}), to name just a few.
They also constitute a~natural continuous time analogue of
random recurrence equations, as studied by de Haan and Karandikar
\cite{HK}.
Lindner and Maller~\cite{LM} have shown that
a generalized Ornstein--Uhlenbeck
process admits a strictly stationary solution which is not degenerate
to a constant process with
a suitable $V_0$ if and only if
%
\begin{equation} \label{eq-convergence}
\int_0^{\infty-} e^{-\xi_{s-}} \, dL_s :=
\lim_{t\to\infty} \int_0^{t} e^{-\xi_{s-}} \, dL_s
\end{equation}
exists and is finite almost surely and not degenerate to
a constant random variable.
The distribution of \textup{(\ref{eq-convergence})} then gives the
unique stationary
distribution.
Here, $\{ (\xi_t, L_t), t \geq0\}$ is another bivariate L\'evy process,
defined in terms of $\{ (\xi_t, \eta_t)\}$ by
\[
L_t = \eta_t + \sum_{0 < s \leq t} \bigl(e^{-(\xi_s - \xi_{s-})}-1\bigr)
(\eta_s
- \eta_{s-})
- t a_{\xi,\eta}^{1,2},
\]
where $a_{\xi,\eta}^{1,2}$ denotes the $(1,2)$-element in the Gaussian
covariance matrix of the L\'evy--Khintchine triplet of $\{(\xi_t,\eta
_t)\}$.
Conversely, $\{(\xi_t,\eta_t)\}$ can be reconstructed from $\{(\xi
_t,L_t)\}$ by
\[
\eta_t = L_t + \sum_{0 < s \leq t} (e^{\xi_s - \xi_{s-}} -1) ( L_s
- L_{s-})
+ t a_{\xi,L}^{1,2}.
\]
Note that, if $\{\xi_t\}$ and $\{\eta_t\}$ are independent, then
$L_t=\eta_t$ for all $t$.
The convergence of integral \textup{(\ref{eq-convergence})}
was characterized by
Erickson and Maller \cite{EM} and generalized by Kondo, Maejima and Sato \cite{KMS}
to the case when $\{(\xi_t,L_t)\}$ is an $\mathbb{R} \times
\mathbb{R}^d$ valued L\'evy process with $d \in\mathbb{N}$.

Suppose now that $\{(\xi_t,L_t)\}$ is a bivariate L\'evy process such that
\textup{(\ref{eq-convergence})} converges almost surely and is
finite, and
denote by
\[
\mu:=\mathcal{L} \biggl(\int_0^{\infty-}
e^{- \xi_{s-} } \, dL_s\biggr)
\]
the distribution of the integral. If $\xi_t = t$ is deterministic,
then it is well known that~$\mu$ is self-decomposable, hence
is infinitely divisible as well as absolutely continuous
(if not degenerate to a Dirac measure, which happens only if
$\{L_t\}$ is also deterministic).
Other cases where $\mu$ is self-decomposable include the case
where
$\{\xi_t\}$ is stochastic, but spectrally negative (cf. Bertoin, Lindner
and Maller \cite{BLM}). On the other hand, as
remarked by Samorodnitsky, $\mu$
is not infinitely
divisible if, for example, $\xi_t = N_t + \alpha t$
with a Poisson process $\{ N_t, t \geq0\}$ and a positive drift
$\alpha> 0$ and $L_t=t$
(cf. Kl{\"u}ppelberg, Lindner and Maller \cite{KLM06}, page 408).
Continuity properties of $\mu$ for general $\{(\xi_t,L_t)\}$ were
studied by Bertoin, Lindner and Maller~\cite{BLM}, who showed
that~$\mu$ cannot have atoms unless $\mu$ is a Dirac measure,
with this degenerate case also being characterized.
Gjessing and Paulsen \cite{GP} derived the distribution of $\mu$ in
a variety of situations; however, in all cases considered
the distribution
turned out to be absolutely continuous.

With these results in mind, it is natural to ask, first, whether
$\mu$ will always be absolutely continuous for general $\{(\xi
_t,L_t)\}$,
unless $\mu$ degenerates to a Dirac measure and, second, what is the
condition for $\mu$ to be infinitely divisible. The
present article will give the negative answer to the first question,
showing many cases of~$\mu$ being continuous-singular and, to
the second question,
provide a necessary and sufficient condition in a restricted class.
Namely we will restrict our attention to the case $(\xi_t,L_t)=
((\log c)N_t,Y_t)$, where $c$ is a constant greater than~1 and
$\{ N_t\}$ and $\{Y_t\}$ are Poisson processes with parameters
$a$ and $b$, respectively, with $\{(N_t,Y_t)\}$ being a bivariate
L\'evy process. Thus we will study in detail
%
\begin{equation} \label{eq-NY}
\mu= \mathcal{L} \biggl(
\int_0^{\infty-} e^{- (\log c) N_{s-}} \, dY_s \biggr) =
\mathcal{L} \biggl(
\int_0^{\infty-} c^{- N_{s-}} \, dY_s \biggr).
\end{equation}
The integral here is
an improper Stieltjes integral pathwise. From the strong law of large
numbers, we see that the integral exists and is finite.
Even in this class the problems of infinite
divisibility and continuity properties turn out to have rich substance.
Let $T$ be the first jump time
of $\{N_t\}$. Then
\[
\int_0^{\infty-} c^{-N_{s-}} \, dY_s = Y_T+
\int_T^{\infty-}c^{-N_{s-}} \, dY_s
= Y_T+c^{-1}\int_0^{\infty-}c^{-N'_{s-}} \, dY'_s,
\]
where $\{(N'_t,Y'_t)\}$ is an independent copy of
$\{(N_t,Y_t)\}$. Hence, letting $\rho=\mathcal L(Y_T)$, we
obtain
%
\begin{equation} \label{eq-SS1}
\widehat{\mu} (z) = \widehat{\rho} (z)   \widehat{\mu} (c^{-1}
z) ,\qquad
z \in\mathbb{R},
\end{equation}
where $\widehat{\mu}(z)$ and $\widehat{\rho}(z)$ denote the
characteristic functions of $\mu$ and $\rho$.
It follows that
\[
\widehat{\mu} (z) = \widehat{\mu} (c^{-k} z)\prod_{n=0}^{k-1}
\widehat{\rho} (c^{-n}z),\qquad   k \in\mathbb{N},
\]
and hence
%
\begin{equation} \label{eq-SS1a}
\widehat{\mu} (z) = \prod_{n=0}^{\infty}\widehat{\rho} (c^{-n} z).
\end{equation}
In general, if a distribution $\mu$ satisfies \textup{(\ref{eq-SS1})}
with some distribution $\rho$,
then~$\mu$ is called \textit{$c^{-1}$-decomposable}.
Our study of the law $\mu$ is based on this\break \mbox{$c^{-1}$-decomposability}.
The expression \textup{(\ref{eq-SS1a})} shows that the law $\rho$ controls
$\mu$. The properties of $c^{-1}$-decomposable distributions
are studied by Wolfe \cite{Wo}, Bunge~\cite{Bu}, Watanabe
\cite{Wa00} and others. In particular, it is known that
any nondegenerate $c^{-1}$-decomposable distribution is
either continuous-singular or absolutely continuous
(Wolfe \cite{Wo}). A distribution $\mu$ is self-decomposable
if and only if~$\mu$ is $c^{-1}$-decomposable for all
$c>1$; in this case $\mu$ and $\rho$ are infinitely divisible.
In general if a distribution $\mu$ satisfies \textup{(\ref{eq-SS1})}
with $\rho$ being infinitely divisible, then~$\mu$ is called
\textit{$c^{-1}$-semi-self-decomposable}.
We note that, when $c=e$ and $\{N_t\}$ and~$\{Y_t\}$ are
independent, Kondo, Maejima and Sato \cite{KMS} recognizes that $\mu$
is $e^{-1}$-decomposable and either continuous-singular or
absolutely continuous.

The L\'evy process $\{(N_t,Y_t)\}$ is a bivariate compound
Poisson process with L\'evy measure concentrated on
the three points $(1,0)$, $(0,1)$ and $(1,1)$ and the
amounts of the measure of these points are denoted by
$u$, $v$ and $w$. Letting $p=u/(u+v+w)$, $q=v/(u+v+w)$
and $r=w/(u+v+w)$ be the normalized L\'evy measure on these
three points, we will see that $\mu$ is determined by
$c$, $q$ and $r$ and~$\rho$ is determined by $q$ and $r$, and hence
denote $\mu=\mu_{c,q,r}$ and $\rho=\rho_{q,r}$. We call
$r$ the \textit{dependence parameter} of $\{(N_t, Y_t) \}$,
since $r=0$ is equivalent to the independence of $\{ N_t \}$
and $\{ Y_t \}$ and $r=1$ means $N_t =Y_t$ for all $t$.
If $r=0$, then $\rho=\mathcal L(Y_T)$ is a geometric distribution,
which is infinitely divisible, and hence $\mu$ is
also infinitely divisible. But, if $r>0$, the situation is
more complicated. In Section~\ref{sec2} we
will give a complete description of the condition of
infinite divisibility of $\mu_{c,q,r}$ and $\rho_{q,r}$ in
terms of their parameters. It will turn out that
infinite divisibility of $\mu_{c,q,r}$ does not depend
on $c$.
It is shown in Niedbalska-Rajba \cite{Ni} that there exists a
$c^{-1}$-decomposable infinitely divisible distribution $\mu$ that
satisfies \textup{(\ref{eq-SS1})} with a noninfinitely-divisible
$\rho$. But, in our case, it will turn out that
$\mu_{c,q,r}$ is infinitely divisible if and only if
$\rho_{q,r}$ is so.
Further, under the condition that $0<q<1$ and $p>0$, it will
turn out that $\rho_{q,r}$ is infinitely divisible if and only
if the dependence parameter is so small that $r\leq pq$.
We also address the problem of infinite divisibility of the
symmetrizations $\mu^{\mathrm{sym}}$ and $\rho^{\mathrm{sym}}$ of
$\mu$ and $\rho$. Infinite divisibility of a
distribution implies that of its symmetrization, but there is
a noninfinitely-divisible distribution whose symmetrization
is infinitely divisible, which is pointed out in pages 81--82
in Gnedenko and Kolmogorov \cite{GK}.
A complete description of infinite divisibility of
$\mu^{\mathrm{sym}}$ and $\rho^{\mathrm{sym}}$ will be given, which
provides new examples of this phenomenon in \cite{GK}.
In the proof of noninfinite-divisibility, we
use three methods: (1) Katti's condition for distributions on
nonnegative integers; (2) L\'evy--Khintchine type representation
of characteristic functions with signed measures in place of
L\'evy measures; (3) representation of the Laplace transforms
of infinitely divisible distributions on $[0,\infty)$ in the
form $e^{-\varphi(\theta)}$ with $\varphi'(\theta)$ being completely
monotone.

Section \ref{sec3} is devoted to the study of continuous-singularity
and absolute continuity of $\mu_{c,q,r}$.
If $q=0$, then it will be shown that
\[
\widehat{\mu_{c,0,r}}(z) = \prod_{n=0}^\infty[ (1-r) +
r e^{ic^{-n} z} ],\qquad
z\in\mathbb{R},
\]
so that $\mu$ is an \textit{infinite Bernoulli convolution}
(usage of this word is not fixed; here we follow
Watanabe \cite{Wa00}).
The question of singularity and absolute continuity of
infinite Bernoulli
convolutions has been investigated by many authors but,
even if $r=1/2$,
characterization of all $c>1$ for which
the distribution is absolutely continuous is an open problem.
See Peres, Schlag and Solomyak \cite{PSS}, Peres and Solomyak \cite{PS},
Watanabe \cite{Wa00} and the references therein.
We shall exclude
the case $q=0$ from our consideration, but we will show that
the notions and
techniques developed in the study of infinite
Bernoulli convolutions and \mbox{$b$-decomposable} measures are effectively
applied.
Here, unlike in the study of infinite divisibility, the parameter
$c$ plays a crucial role. If $c$ has an algebraic property of
being a Pisot--Vijayaraghavan (P.V.) number, then
we will show that $\mu_{c,q,r}$
is continuous-singular under the condition that
$q>0$. For example, all integers greater than one
and some irrationals such as $(1+\sqrt{5})/2$ are P.V. numbers.
On the other hand, if $c$ is the reciprocal of
a Peres--Solomyak (P.S.) number, then it will be shown that
there are positive constants $C_1$ and $C_2$ such that $\mu_{c,q,r}$
is absolutely continuous with bounded continuous density
whenever $q\geq C_1p \geq C_2 r$. It is known that
Lebesgue almost all reals in $(1,\infty)$ are reciprocals
of P.S. numbers.
In general, under the condition $0<q<1$, we can estimate
$\dim (\mu_{c,q,r})$, the
Hausdorff dimension of $\mu_{c,q,r}$ defined as
the infimum of the Hausdorff dimensions of $E$
over all Borel sets $E$ satisfying $\mu_{c,q,r} (E) =1$
(in some papers, including \cite{Wa00}, this is called upper
Hausdorff dimension and denoted by $\dim^*$).
Using a powerful theorem of Watanabe \cite{Wa00}
for any $c^{-1}$-decomposable distribution satisfying
\textup{(\ref{eq-SS1})} with a discrete distribution $\rho$,
we see that $\dim (\mu_{c,q,r})\leq H(\rho_{q,r})/\log c$,
where $H(\rho_{q,r})$ is the
entropy of $\rho_{q,r}$. It follows that $\mu_{c,q,r}$
is continuous-singular, if $H(\rho_{q,r})/\log c<1$.
Thus, for any $c>1$, there is a positive constant $C_3$
such that $\mu_{c,q,r}$
is continuous-singular whenever $q>0$ and
$\max\{q,r\}\leq C_3 p$.

In Section \ref{sec3} we also study, in the case where $\mu_{c,q,r}$
is infinitely divisible, continuity properties
of the convolution power
$(\mu_{c,q,r})^{t*}$ of $\mu_{c,q,r}$, that is, the
distribution at time $t$ of the L\'evy process
associated with $\mu_{c,q,r}$. It is shown that
if $c$ is a P.V. number, then $(\mu_{c,q,r})^{t*}$ is
continuous-singular for all $t>0$, while, if $c$ is the
reciprocal of a P.S. number, then there is $t_1=t_1(c,q,r)
\in(0,\infty)$ such that $(\mu_{c,q,r})^{t*}$ is
continuous-singular for all $t\in(0,t_1)$ and absolutely
continuous for all $t\in(t_1,\infty)$.
Thus the present paper provides a new class
of L\'evy processes with a remarkable time evolution in
distribution. See Section 27 in Sato \cite{Sa} and Watanabe's
survey \cite{Wa01} for such time evolution.
We emphasize that here the distribution $\mu_{c,q,r}$
arises naturally as the stationary distribution of a
generalized Ornstein--Uhlenbeck process.

The case of $\{ N_t \}$ and $\{ Y_t \}$
being independent (i.e., $r=0$) is of special interest.
The properties of $\mu_{c,q,0}$ are included in the
results of Section \ref{sec3} mentioned above. As explicit
examples, $\mu_{e,q,0}$ with $c=e$ is continuous-singular
if $q\leq1-\log2\approx0.30685$; $\mu_{c,1/2,0}$
with $q=1/2$ is continuous-singular if $c>4$.
We can prove more results for $\mu_{c,q,0}$ than for
general $\mu_{c,q,r}$, since the L\'evy measure of
$\mu_{c,q,0}$ is increasing with respect to $q$.
Thus, for any $c>1$, there exists $q_1$ with
$0<q_1\leq1$ such that $\mu_{c,q,0}$ is continuous-singular
for all $q\in(0,q_1)$ and absolutely continuous for all
$q\in(q_1,1)$. It will be shown that $q_1=1$ for any P.V. number $c>1$
and that $q_1<1$ whenever $c$ is the reciprocal of a P.S. number, so
that $q_1<1$ for Lebesgue almost all $c>1$.


Throughout the paper, the set of all
positive integers will be denoted
by $\mathbb{N} = \{1,2 ,3, \ldots\}$, while we set $\mathbb{N}_0
= \mathbb{N} \cup\{ 0 \}$.
The set of integers is denoted by $\mathbb{Z}$.
The Dirac measure at a point $x$ will be denoted by $\delta_x$.
For general definitions and properties
regarding L\'evy processes and infinitely divisible distributions,
we refer to Sato \cite{Sa}.

\section{Necessary and sufficient conditions for infinite divisibility}
\label{sec2}

Suppose that $\{ (N_t, Y_t), t\geq0\}$
is a bivariate L\'evy process such that $\{N_t\}$ is a Poisson process
with parameter $a> 0$ and $\{Y_t\}$ is a Poisson process with parameter
$b > 0$.
It then follows easily that $\{(N_t,Y_t)\}$ has no Gaussian part,
no drift, and a L\'evy measure $\nu_{(N,Y)}$ concentrated on
the set $\{ (1,0), (0,1), (1,1) \}$, consisting of three points
(e.g., \cite{Sa}, Proposition 11.10).
Denote
\[
u := \nu_{(N,Y)}( \{ (1,0) \}),\qquad
v := \nu_{(N,Y)} ( \{ (0,1) \}), \qquad
w := \nu_{(N,Y)} ( \{ (1,1) \} ).
\]
Then $u,v,w \geq0$,
$u+w = a$ and $v+w=b$. Let
\[
p:=\frac{u}{u+v+w},\qquad  q:=\frac{v}{u+v+w},\qquad  r:=\frac{w}{u+v+w},
\]
so that $p,q,r\in[0,1]$, $p+q+r=1$, $p+r>0$ and $q+r>0$.
These give the normalized L\'evy measure on the
three points.
The two processes $\{N_t\}$ and $\{Y_t\}$ are independent if and only if
$r=0$.
If $r=1$, then $N_t=Y_t$ for all $t$ with probability one.
So we call $r$ the \textit{dependence parameter} of
$\{(N_t,Y_t)\}$.
The law $\mu$ in \textup{(\ref{eq-NY})} depends on $c$, $u$, $v$ and
$w$. But it
will turn out (Proposition \ref{prop-3-decom}) that $\mu$ depends only
on $c$, $q$ and $r$. Thus, for $c>1$ denote
%
\begin{equation} \label{eq-mu2}
\mu_{c,q,r} := \mathcal{L} \biggl( \int_0^{\infty-} c^{-N_{s-}} \,
dY_s \biggr).
\end{equation}
If $r=1$, then
\[
\int_0^{\infty-} c^{-N_{s-}} \, dY_s=
\int_0^{\infty-} c^{-N_{s-}} \, dN_s =
\sum_{j=0}^\infty c^{-j} = \frac{c}{c-1},
\]
which is degenerate to a constant.
So, from now on, we assume that $p+q>0$
in addition to the above-mentioned conditions
$p+r > 0$ and $q+r>0$. That is, $p,q,r<1$. In this section
we are interested in whether $\mu_{c,q,r}$ is infinitely
divisible or not. It is also of interest whether the
symmetrization $(\mu_{c,q,r})^{\mathrm{sym}}$ of $\mu_{c,q,r}$ is
infinitely divisible or not.
Recall that the symmetrization $\mu^{\mathrm{sym}}$ of a
distribution~$\mu$ is defined to be the distribution with characteristic function
$|\widehat\mu(z)|^2$.
Infinite divisibility of $\mu$ implies that
of $\mu^{\mathrm{sym}}$, but the converse is not true, as is
mentioned in the \hyperref[sec1]{Introduction}.

We define $\rho_{q,r}$ in the following
way: If $q > 0$, denote by
$\sigma_q$ a geometric distribution with parameter $1-q$, that is,
$\sigma_q (\{ k \}) = (1-q) q^k$ for $k=0,1,\ldots$ and denote
%
\begin{equation} \label{eq-rho-1}
\rho= \rho_{q,r} := (1 + r/q) \sigma_q - (r/q)\delta_0,
\end{equation}
so that $\rho_{q,r}$ is a probability distribution concentrated
on $\mathbb{N}_0$ with
%
\begin{equation} \label{eq-rho-1a}
\rho_{q,r}(\{0\})=(1 + r/q)(1-q) - (r/q)=p ;
\end{equation}
if $q = 0$, let $\rho_{0,r}$ be a Bernoulli
distribution with parameter $r\in(0,1)$, that is,
%
\begin{equation} \label{eq-rho-1b}
\rho_{0,r}(\{1\})=1-\rho_{0,r}(\{0\})=r.
\end{equation}

\begin{proposition} \label{prop-3-decom}
We have
%
\begin{equation} \label{eq-SS3}
\widehat{\mu}_{c,q,r} (z) = \widehat{\rho}_{q,r} (z)
\widehat{\mu}_{c,q,r} (c^{-1} z), \qquad  z\in\mathbb{R}.
\end{equation}
In particular, $\mu_{c,q,r}$ is $c^{-1}$-decomposable and
determined by $c$, $q$ and $r$.
\end{proposition}

\begin{pf}
As is explained in the \hyperref[sec1]{Introduction}, we have only to
show that $\mathcal L(Y_T)=\rho_{q,r}$, where $T$ is
the time of the first jump of $\{N_t\}$, that is, the time of the
first jump of $\{(N_t,Y_t)\}$ with size in $\{ (1,0), (1,1) \}$.
Let $S_i$ be the size of the $i$th jump of
$\{(N_t,Y_t)\}$. Then we have for $k \geq1$
\begin{eqnarray*}
Y_{T} = k\quad &\Longleftrightarrow\quad &[ S_1 = \cdots=
S_{k-1} = ( 0, 1 ),
S_k = ( 1 , 1) ]\\
&&
\mbox{or}\\
&&[ S_1 = \cdots= S_{k} = (0 , 1),
S_{k+1} = ( 1 , 0) ],
\end{eqnarray*}
as well as
\[
Y_{T} = 0 \quad\Longleftrightarrow\quad  S_1 = (1 , 0 ).
\]
Since
\[
P[ S_i = (1 , 0 ) ] = p,\qquad
P[ S_i = (0 , 1) ] = q, \qquad  P[S_i = ( 1 , 1) ] = r,
\]
it follows that $P(Y_{T} = 0) = p$
and, for $k \geq1$, $P( Y_{T} = k) = q^{k-1}r+q^k p$.
From this it follows easily that $\mathcal L(Y_T)=\rho_{q,r}$ for $q>0$,
while it
is a Bernoulli distribution with parameter $r$ for $q=0$.
\end{pf}

We can now formulate criteria when $\rho_{q,r}$ and $\mu_{c,q,r}$
and their symmetrizations
are infinitely divisible. As is seen in \textup{(\ref{eq-SS1a})}, infinite
divisibility of $\rho_{q,r}$ implies that of $\mu_{c,q,r}$.
Similarly, infinite
divisibility of $(\rho_{q,r})^{\mathrm{sym}}$ implies that of
$(\mu_{c,q,r})^{\mathrm{sym}}$. The converse of these two implications
is by no means clear, as we know Niedbalska-Rajba's example
mentioned in the \hyperref[sec1]{Introduction}. However the following theorem will
say that the converse is true for $\mu_{c,q,r}$ and $\rho_{q,r}$ and
for $(\mu_{c,q,r})^{\mathrm{sym}}$ and $(\rho_{q,r})^{\mathrm{sym}}$.
Thus infinite divisibility of $\mu_{c,q,r}$ does not depend on $c$.
Another remarkable consequence is that $(\mu_{c,q,r})^{\mathrm{sym}}$
can be infinitely divisible without $\mu_{c,q,r}$ being
infinitely divisible and that $(\rho_{q,r})^{\mathrm{sym}}$ can be
infinitely divisible without $\rho_{q,r}$ being
infinitely divisible.

\begin{theorem} \label{thm-id}
Assume that the parameters $p,q,r$ satisfy $p,q,r<1$. Let $c>1$.
For $\mu_{c,q,r}$ defined in \textup{(\ref{eq-mu2})}
and $\rho_{q,r}$ in
\textup{(\ref{eq-rho-1})} and \textup{(\ref{eq-rho-1b})}, the
following hold true:\vspace*{-6pt}

\begin{longlist}[(a)]
\item[(a)] If $p=0$, then $\rho_{q,r}$ and $\mu_{c,q,r}$
are infinitely divisible.
\item[(b)] If $p > 0$ and $q> 0$, then the following conditions
are equivalent:\vspace*{-6pt}
\begin{longlist}[(a)\,\ (iii)]
\item $\mu_{c,q,r}$ is infinitely divisible.
\item $\rho_{q,r}$ is infinitely divisible.
\item $r\leq pq$.
\end{longlist}\vspace*{-6pt}
\item[(c)] If
$p>0$, $q> 0$ and $r>pq$, then the
following conditions are equivalent:\vspace*{-6pt}
\begin{longlist}[(a)\ \,(iii)]
\item $(\mu_{c,q,r})^{\mathrm{sym}}$ is infinitely divisible.
\item  $(\rho_{q,r})^{\mathrm{sym}}$ is infinitely divisible.
\item  $p\leq qr$.
\end{longlist}\vspace*{-6pt}
\item[(d)] If $q=0$, then none of $\rho_{q,r}$, $\mu_{c,q,r}$,
$(\rho_{q,r})^{\mathrm{sym}}$ and $(\mu_{c,q,r})^{\mathrm{sym}}$
is infinitely divisible.
\end{longlist}
\end{theorem}

In the proof, we will first settle the
question of infinite divisibility of $\rho_{q,r}$ and
$(\rho_{q,r})^{\mathrm{sym}}$.

\begin{lemma} \label{lem-rho-id}
Assume $q>0$ and let $\rho=\rho_{q,r}$.
Then the following hold true:

\textup{(a)} If $r\leq pq$, or if $p=0$,
then $\rho$ is infinitely divisible.

\textup{(b)} If $r>pq$ and $p>0$, then $\rho$ is not infinitely divisible.

\textup{(c)} Assume that $r>pq$ and $p>0$. Then
$\rho^{\mathrm{sym}}$ is infinitely divisible if and only if
%
\begin{equation} \label{eq-rho-sym}
p\leq qr.
\end{equation}
\end{lemma}

We remark that if $0 \leq\alpha\leq1$, then
$(1-\alpha) \sigma_q + \alpha\delta_0$
is infinitely divisible, since convex combinations of two
geometric distributions are infinitely divisible (see pages 379--380
in Steutel and van Harn \cite{SvH}), and the Dirac measure $\delta_0$
is a limit of geometric distributions. Assertions \textup{(a)} and
\textup{(b)}
show to what extent this fact can be generalized to negative $\alpha$.

\begin{pf*}{Proof of Lemma \protect\ref{lem-rho-id}}
Since $\widehat\sigma_q(z)=(1-q)/(1-qe^{iz})$, we have
%
\begin{equation} \label{eq-cf}
\widehat\rho(z)=\frac{p+re^{iz}}{1-qe^{iz}},\qquad   z\in\mathbb{R}.
\end{equation}

\textup{(a)} If $p=0$, then $\rho(\{0\})=0$ and $\rho(\{k\})
=(1-q) q^{k-1}$ for $k=1,2,\ldots$ and thus~$\rho$ is a geometric
distribution
translated by $1$, hence infinitely divisible. So assume that
$r\leq pq$.
Then $p>0$, recalling that $p+r>0$.
Since $p = (1-q)/(1+r/p)$,
it follows from \textup{(\ref{eq-cf})} that
\[
\widehat{\rho} (z)
=\exp\biggl[\log(1-q) - \log\biggl(1+\frac{r}{p}
\biggr) +\log\biggl(1+\frac{r}{p} e^{iz}\biggr)-\log(1-qe^{iz})
\biggr].
\]
Hence
%
\begin{equation}\label{eq-cf2}
\widehat\rho(z)=\exp\Biggl[\sum_{k=1}^{\infty}(e^{ikz}-1)\frac{q^k}{k}
\biggl(1-\biggl(-\frac{r}{pq}\biggr)^k\biggr)\Biggr].
\end{equation}
Recall that $r/(pq)\leq1$. It follows that
$\rho$ is infinitely divisible with L\'evy measure
$\nu_{\rho}(\{k\})=k^{-1} q^k (1-(-r/(pq))^k)$, $k=1,2,\ldots,$ and
drift $0$.

\textup{(b)} Now assume that $r>pq$ and $p>0$. By Katti's criterion
(\cite{K67} or Corollary~51.2 of \cite{Sa}), a
distribution $\sum_{n=0}^{\infty} p_n \delta_n$ with $p_0>0$ is
infinitely divisible if and only if there are $q_n\geq0$,
$n=1,2,\ldots
,$ such that
\[
np_n=\sum_{k=1}^n kq_k p_{n-k},\qquad n=1,2,\ldots .
\]
In fact, the equations above determine $q_n$, $n=1,2,\ldots,$
successively in a
unique way; infinite divisibility of $\sum_{n=0}^{\infty} p_n \delta
_n$ is
equivalent
to nonnegativity of all $q_n$. Now let $p_n=\rho(\{n\})$.
The first two equations are
$p_1=q_1 p_0$ and $2p_2=q_1 p_1+2q_2 p_0$.
Hence $q_1=p_1/p_0 > 0$, but
\[
q_2=\frac{2p_2-q_1 p_1}{2p_0}
=\frac{(1+r/q)(1-q) q^2}{2p^2} [1-q-(r/q)(1+q)]
<0,
\]
since $r>pq$. This shows that $\rho$ is not infinitely divisible.

\textup{(c)} Assume again that $r>pq$ and $p>0$. From \textup{(\ref{eq-cf})}
it can be seen that $\widehat{\rho}$ will have a real zero if and
only if $p=r$.
In that case, $|\widehat{\rho}|^2$ will also
have a real zero, and hence $\rho^{\mathrm{sym}}$ cannot be infinite divisible
in agreement with the fact that \textup{(\ref{eq-rho-sym})} is
violated for $p=r$.
So in the following we assume that $p\neq r$.
From \textup{(\ref{eq-cf})} we have
\[
\log(|\widehat\rho(z)|^2)=\log(p^2+2pr\cos z+r^2)-\log(1-2q\cos z+q^2).
\]
Write
\[
A=\frac{2pr}{p^2+r^2},\qquad
B=\frac{2q}{1+q^2},\qquad
C=\frac{p^2+r^2}{1+q^2}.
\]
Then $0<A<1$, $0<B<1$, and $C>0$ (recall that $0<q<1$ and $p\neq r$)
and we obtain
\begin{eqnarray*}
\log(|\widehat\rho(z)|^2)&=&\log C + \log(1 + A \cos z) - \log(1- B
\cos
z)\\
&=&\log C-\sum_{k=1}^{\infty}k^{-1}(-A)^k
\cos^k z+\sum_{k=1}^{\infty}k^{-1}B^k\cos^k z\\
&=&\log C+\sum_{k=1}^{\infty}k^{-1}2^{-k}\bigl(-(-A)^k+B^k\bigr)
\sum_{l=0}^k \pmatrix{k\cr l}\cos(k-2l)z,
\end{eqnarray*}
since
$
\cos^k z=2^{-k}\sum_{l=0}^k {k\choose l}\cos(k-2l)z.
$
Letting $z=0$, we represent $\log C$ by $A$ and $B$ and get
\[
\log(|\widehat\rho(z)|^2)=\sum_{k=1}^{\infty}D_k
\sum_{l=0}^k \pmatrix{k\cr l}\bigl(\cos(k-2l)z-1\bigr),
\]
where
%
\begin{equation} \label{eq-Dk}
D_k=k^{-1}2^{-k}\bigl(-(-A)^k+B^k\bigr).
\end{equation}
Then we get, with $\lfloor(k-1)/2\rfloor$ denoting
the largest integer not exceeding $(k-1)/2$,
\begin{eqnarray*}
\log(|\widehat\rho(z)|^2)&=&2\sum_{k=1}^{\infty} D_k
\sum_{l=0}^{\lfloor(k-1)/2\rfloor} \pmatrix{k\cr l}\bigl(\cos
(k-2l)z-1\bigr)\\
&=&2\sum_{k=1}^{\infty} D_k \sum_m \pmatrix{k\cr(k-m)/2}
(\cos mz-1),
\end{eqnarray*}
where $m$ runs over $k, k-2, \ldots,3, 1$, if $k$ is odd $\geq1$
and over $k, k-2, \ldots,4, 2$, if~$k$ is even $\geq2$.
Since
$
\sum_{k=1}^{\infty} 2^k |D_k|\leq\sum_{k=1}^{\infty} k^{-1}
(A^k+B^k)<\infty
$,
we can change the order of summation and obtain
%
\begin{equation}\label{eq-1.3}
\log(|\widehat\rho(z)|^2)=2\sum_{m=1}^{\infty} E_m(\cos mz-1)
\end{equation}
with
%
\begin{equation} \label{eq-Em}
E_m = \sum_{h=0}^\infty D_{m+2h} \pmatrix{m+2h\cr h}.
\end{equation}
This means that
%
\begin{equation}\label{eq-1.4}
\log(|\widehat\rho(z)|^2)=\int_{\mathbb{R}} \bigl(e^{ixz}-1-ixz
1_{(-1,1)}(x)\bigr) \nu(dx),
\end{equation}
where $\nu$ is the symmetric signed measure
%
\begin{equation}\label{eq-1.5}
\nu=\sum_{m=1}^{\infty} E_m(\delta_m+\delta_{-m}).
\end{equation}
Let $F=r/(pq)$. Then $F>1$. A simple calculation
then shows that $A \leq B$
if and only if $F-1\leq q^2(F^2-F)$,
which is equivalent to $1\leq q^2 F$, that is,
\textup{(\ref{eq-rho-sym})}.
Now, if~\textup{(\ref{eq-rho-sym})} holds, then $A\leq B$ and hence
$D_k\geq0$ for all $k$, which implies $E_m\geq0$ for all~$m$
and $\rho^{\mathrm{sym}}$ is infinitely divisible with
the L\'evy--Khintchine representation~\textup{(\ref{eq-1.4})} with
\textup{(\ref{eq-1.5})}.
If \textup{(\ref{eq-rho-sym})} does not hold, then $A> B$,
$D_k<0$ for all even $k$, and $E_m<0$ for all even $m$, which
implies, by \textup{(\ref{eq-1.4})} and \textup{(\ref{eq-1.5})},
that $\rho^{\mathrm
{sym}}$ is
not infinitely divisible (see Exercise 12.3 of \cite{Sa}).
\end{pf*}

\begin{pf*}{Proof of Theorem \protect\ref{thm-id}}
Write $\mu=\mu_{c,q,r}$ and $\rho=\rho_{q,r}$.
(a) Suppose $p=0$. Then $\rho$ is infinitely
divisible by Lemma \ref{lem-rho-id}, and hence so is $\mu$ by
\textup{(\ref{eq-SS3})}.

(b) Suppose that $p,q > 0$. Under these conditions,
the equivalence of (ii) and~(iii)
follows from Lemma \ref{lem-rho-id}. Further, (ii) implies (i) by
\textup{(\ref{eq-SS3})}, so that it remains to show that (i) implies (iii).
For that, suppose that $r>pq$, and in order to show
that $\mu$ is not infinitely divisible,
we will distinguish three cases: $p=r$, $p>r$ and
$p<r$.
The first case is easy, but in the second and third cases we
have to use rather involved arguments resorting to different
conditions that guarantee noninfinite-divisibility.

\textit{Case} 1: Suppose that $p=r$. Then
$\widehat{\rho}$ will have a real zero as argued in the proof
of Lemma \ref{lem-rho-id}(c). By \textup{(\ref{eq-SS3})},
$\widehat{\mu}$ will also
have a real zero, so that $\mu$ cannot be infinitely divisible.

\textit{Case} 2: Suppose that $p>r$. Then
$\widehat{\rho}$ can be expressed as in \textup{(\ref{eq-cf2})}
with the same
derivation.
Together with \textup{(\ref{eq-SS3})} and \textup{(\ref{eq-SS1a})}
this implies
%
\begin{equation} \label{eq-inf-prod}\qquad
\widehat\mu(z)=
\exp\Biggl[\sum_{n=0}^{\infty}
\sum_{m=1}^{\infty}(e^{imc^{-n}z}-1)\frac
{1}{m}\bigl(q^{m}-(-r/p)^{m}\bigr)\Biggr],
\qquad  z\in\mathbb{R}.
\end{equation}
Absolute convergence of this double series follows from $c>1$,
$q < 1$ and $r/p<1$.
Define the real numbers $a_m$, $m\in\mathbb{N}$, and the signed measure
$\nu$ by
\[
a_m := \frac{1}{m} \bigl(q^m - (-r/p)^m\bigr)\quad  \mbox{and}\quad
\nu:= \sum_{n=0}^\infty\sum_{m=1}^\infty a_m   \delta_{c^{-n}m}.
\]
It follows that $\widehat{\mu}$ in \textup{(\ref{eq-inf-prod})}
has the same form as the L\'evy--Khintchine
representation with the signed measure $\nu$ in place of a L\'evy
measure, so that
infinite divisibility of $\mu$ is equivalent to the signed measure
$\nu$
having negative part $0$; see Exercise 12.3 in \cite{Sa}.
Thus, to show that $\mu$ is not infinitely divisible, we will show that
there is a point $x$ such that $\nu(\{ x \}) < 0$.
Since $r/p> q$,
it follows that $a_m < 0$ if $m$ is even and that $a_m > 0$ if $m$ is odd.
If $c^k$ is irrational for all $k\in\mathbb{N}$, then the points
$c^{-n}m$ with
$n\in\mathbb{N}_0$ and $m\in\mathbb{N}$ are distinct, which implies
$\nu(\{
c^{-n}m\})<0$
for all even $m$ and $\mu$ is not infinitely divisible.

Suppose that $c^k$ is rational for some $k\in\mathbb{N}$. Let $k_0$
be the
smallest such $k$ and write $c^{k_0} = \alpha/\beta$ with $\alpha
,\beta\in
\mathbb{N}$
such that $\alpha$ and $\beta$ have no common divisor. Let
$f$ be the largest $t\in\mathbb{N}_0$ such that $2^t$
divides $\beta$. Let $m$ be even. Denote
\begin{eqnarray*}
G_m & := & \{ (n',m') \in\mathbb{N}_0 \times\mathbb{N}\dvtx
c^{-n'} m' = m,
\  m'  \mbox{ odd}\},\\
H_m & := & \{ (n',m') \in\mathbb{N}_0 \times\mathbb{N}\dvtx
c^{-n'} m' = m,
\  m'  \mbox{ even}\}.
\end{eqnarray*}
Then
%
\begin{equation} \label{eq-E2.1}
\nu( \{ m \} ) = \sum_{(n',m') \in G_m \cup H_m} a_{m'}   \leq
  a_m + \sum_{(n',m') \in G_m} a_{m'}.
\end{equation}
We claim that the set $G_m$ contains at most one element. To show
this, let $(n',m')\in G_m$. Then $n'\neq0$ and $c^{n'}=m'/m$ and thus
$n'=lk_0$ for some $l\in\mathbb{N}$. Then $m'/m=(\alpha/\beta)^l$
and hence
$\beta^l$ divides $m$ and $m/\beta^l$ is odd. Thus $m=2^{lf}m''$ with
some odd integer $m''$. It follows that $lf$ is determined by $m$.
Hence $l$ is determined by $m$ and $c$. Hence $n'$
is determined by $m$ and $c$, which shows that
$G_m$ contains at most one element. It also follows that $f\geq1$
whenever $G_m\neq\varnothing$ for some even $m$.
If there is some even $m$ such that $G_m = \varnothing$, then $\nu(\{ m
\})
< 0$ by \textup{(\ref{eq-E2.1})}, and we are done. So suppose from
now on that
$G_m \not= \varnothing$ for every even $m$.
Let $m_j = 2^{j f}$ for $j=1,2,\ldots.$ The argument above shows that
the unique element $(n_j',m_j')$ in $G_{m_j}$ is given by $n_j' = j
k_0$ and
$m_j'=c^{j k_0} m_j$.
Noting that $0<q<r/p<1$ and $c>1$, choose $j$ so large that $q^{m_j}\le2^{-1}
(r/p)^{m_j}$ and $m'_j=m_j c^{jk_0}>2m_j$. Then
\[
a_{m_j}=\frac{1}{m_j}\bigl(q^{m_j}-(r/p)^{m_j}\bigr)\le-\frac{1}{2m_j}(r/p)^{m_j}
\]
and
\[
a_{m'_j}=\frac{1}{m'_j}\bigl(q^{m'_j}+(r/p)^{m'_j}\bigr)\le\frac{3}{2m'_j}(r/p)^{m'_j}
<\frac{3}{4m_j}(r/p)^{2m_j}.
\]
Thus
\[
\nu(\{m_j\})\le a_{m_j}+a_{m'_j}\le\frac{1}{2m_j}(r/p)^{m_j}
\bigl(-1+(3/2)(r/p)^{m_j}\bigr)<0
\]
for large enough $j$, showing that $\mu$ is not infinitely divisible
under the conditions of Case 2.

\textit{Case} 3: Suppose that $p<r$ and, by way of
contradiction, assume that
$\mu$ is infinitely divisible.
Denote by $L_{\mu}(\theta)= \int_{\mathbb{R}}e^{-\theta x} \mu(dx)$,
$\theta\geq0$, the Laplace transform of $\mu$.
Then $L_{\mu}(\theta)=e^{-\varphi(\theta)}$
where $\varphi$ has a completely monotone derivative $\psi(\theta)$ on
$(0,\infty)$, that is, $(-1)^n \psi^{(n)}(\theta)\geq0$ on
$(0,\infty)$
for $n=0,1,\ldots$ (see Feller \cite{Fe}, page 450).
By \textup{(\ref{eq-SS3})} and \textup{(\ref{eq-cf})} we have
%
\begin{equation} \label{eq-11-1}
\varphi(\theta)=-\log L_{\mu}(\theta)=-\sum_{n=0}^{\infty}\log
\frac
{p+rf_n(\theta)}{
1-qf_n(\theta)},
\end{equation}
where $f_0(\theta)=e^{-\theta}$ and
$f_n(\theta)=\exp(-c^{-n}\theta)=f_0(c^{-n}\theta)$,
$n=1,2,\ldots.$ Convergence of the summation in
\textup{(\ref{eq-11-1})} is easily established.
Since $\psi= \frac{d}{d\theta} \varphi$ is completely monotone,
so is $\theta\mapsto c^{-1} \psi(c^{-1} \theta) = \frac{d}{d\theta}
(\varphi(c^{-1} \theta))$. Consider the function
\[
\xi(\theta) =
\frac{p}{p + r e^{-\theta}},\qquad \theta\in(0,\infty).
\]
Then $1/(1-qf_0(\theta))-\xi(\theta)$ is the difference of two
completely monotone functions, because
\[
\frac{d}{d\theta} \bigl(\varphi(\theta) - \varphi(c^{-1} \theta) \bigr)
=\frac{d}{d\theta}
\biggl( - \log\frac{p+ r f_0 (\theta)}{1-q f_0(\theta)}\biggr)
=\frac{1}{1-q f_0(\theta)} -\xi(\theta).
\]
Since $1/(1-qf_0(\theta))=\sum_{k=0}^{\infty}q^k e^{-k\theta}$ is
completely monotone, $\xi(\theta)$ is itself the difference of two
completely monotone functions.
Applying Bernstein's theorem, there must exist a signed measure
$\sigma$ on $[0,\infty)$ such that $\int_{[0,\infty)} e^{-\theta x}
|\sigma|(dx)<\infty$ and
$\xi(\theta) = \int_{[0,\infty)} e^{-\theta x} \sigma(dx)$
for all $\theta\in(0,\infty)$. However,
introducing the signed measure
$\tau:= \sum_{k=0}^\infty(-r/p)^k \delta_k$,
we have
\[
\xi(\theta) = \sum_{k=0}^\infty\biggl( -\frac{r}{p} e^{-\theta}
\biggr)^k = \int_{[0,\infty)} e^{-\theta x} \tau(dx),
\]
if $\theta>\theta_0:=\log(r/p)$. Thus $e^{-\theta_0 x}\sigma(dx)$ and
$e^{-\theta_0 x}\tau(dx)$ have a common Laplace transform
$\xi(\theta_0+\theta)$, $\theta>0$. Now from the
uniqueness theorem in Laplace transform theory (page 430 of
Feller \cite{Fe}) combined with the Hahn--Jordan decomposition of
signed measures, it follows that $e^{-\theta_0 x}\sigma(dx)
=e^{-\theta_0 x}\tau(dx)$, that is,
$\sigma= \tau$.
But $\int_{[0,\infty)}e^{-\theta x} |\tau|(dx)=\infty$ for
$0<\theta\leq\theta_0$, contradicting
the property of $\sigma$. This finishes the proof of (b).

(c) Suppose that $p,q> 0$ and that $r>pq$.
The equivalence of (ii) and (iii)
then follows from Lemma \ref{lem-rho-id}, and (ii) implies (i) by
\textup{(\ref{eq-SS3})}, so that it remains to show that (i) implies (iii).
If $p=r$, then $|\widehat{\rho}|^2$ and hence $|\widehat{\mu}|^2$
have real zeros
as shown in the proof of Lemma \ref{lem-rho-id}(c) and
$\mu^{\mathrm{sym}}$ is not infinitely divisible.
Hence we can assume that $p \neq r$.
With $A, B$, $D_k$ and $E_m$ as in the proof of Lemma \ref
{lem-rho-id}(c),
it follows from $|\widehat{\mu}(z)|^2 = \prod_{n=0}^\infty|\widehat
{\rho
}(c^{-n}z)|^2$
and \textup{(\ref{eq-1.3})} that
%
\begin{equation}\label{eq-1.13}
\log(|\widehat\mu(z)|^2)=2\sum_{n=0}^{\infty}\sum_{m=1}^{\infty} E_m
\bigl(\cos(mc^{-n}z)-1\bigr).
\end{equation}
Since
\begin{eqnarray*}
&&2\sum_{n=0}^{\infty}\sum_{m=1}^{\infty}|E_m| |\cos(mc^{-n}z)-1|\\
&&\qquad =\sum_{n=0}^{\infty}\sum_{k=1}^{\infty}|D_k|\sum_{l=0}^k
\pmatrix{k\cr l}\bigl|\cos\bigl((k-2l)c^{-n}z\bigr)-1\bigr|\\
&&\qquad \le\sum_{n=0}^{\infty}\sum_{k=1}^{\infty}|D_k| 2^k (kc^{-n}z)^2
\le z^2\sum_{n=0}^{\infty}c^{-2n} \sum_{k=1}^{\infty} k(A^k+B^k)
<\infty,
\end{eqnarray*}
we can consider the right-hand side of \textup{(\ref{eq-1.13})} as an integral
with respect to a signed measure. Thus
%
\begin{equation}\label{14}
\log(|\widehat\mu(z)|^2)=\int_{\mathbb{R}}\bigl(e^{ixz}-1-ixz
1_{(-1,1)}(x)\bigr) \widetilde
\nu(dx),
\end{equation}
where $\widetilde\nu$ is the symmetric signed measure
%
\begin{equation}\label{15}
\widetilde\nu=\sum_{n=0}^{\infty}\sum_{m=1}^{\infty} E_m(\delta_{mc^{-n}}
+\delta_{-mc^{-n}}).
\end{equation}
Now suppose that $p>qr$.
As observed in the proof of Lemma \ref{lem-rho-id}(c), this is equivalent
to $A > B$.
In order to show that $\mu$ is not infinitely divisible, we use
Exercise~12.3 of \cite{Sa} again. We need to show that
$\widetilde{\nu}$ has a nontrivial negative part.
Recall that $E_m > 0$ for all odd $m$
and $E_m < 0$ for all even $m$.
If $c^k$ is irrational for all $k\in\mathbb{N}$, then $\widetilde\nu
(\{m\}
)=E_m<0$ for
even $m$. Hence, suppose that $c^k$ is rational for some $k\in\mathbb{N}$.
We first estimate $E_m$. Since
${m+2h\choose h} \leq2^{m+2h}$,
it follows from
\textup{(\ref{eq-Dk})} and \textup{(\ref{eq-Em})} that
%
\begin{equation} \label{eq-Em-odd}
|E_m| \leq\sum_{h=0}^\infty\frac{1}{m+2h} 2 A^{m+2h} \leq
\frac{2A^m}{m(1-A^2)}, \qquad  m\in\mathbb{N}.
\end{equation}
Choose $\gamma\in(0,1)$ such that $A/\gamma< 1$, and choose
$\alpha\in\mathbb{N}$ such that $(\alpha+1/2)/(\alpha+1) \geq
\gamma$. By
Stirling's formula,
there exists a constant $d_1 > 0$ such that for every $m\in\mathbb{N}$,
\begin{eqnarray*}
\pmatrix{m+2\alpha m\cr\alpha m} & \geq & d_1 \biggl( \frac{m+2\alpha
m}{(m+\alpha
m)\alpha m}
\biggr)^{1/2}
\frac{ (m+2\alpha m)^{m+2\alpha m}}{ (m+\alpha m)^{m+\alpha m} (\alpha
m)^{\alpha m}} \\
& \geq & \frac{d_1}{(\alpha m)^{1/2}}   (2\gamma)^{m+2\alpha m}.
\end{eqnarray*}
Since $D_k < 0$ for every even $k$, we conclude
\begin{eqnarray}\label{eq-Em-even}
|E_m| & \geq& |D_{m+2\alpha m}| \pmatrix{m+2\alpha m\cr\alpha m}
\nonumber\\[-8pt]\\[-8pt]
& \geq& \frac{d_1}{ (\alpha m)^{1/2} (m+2\alpha m)}
(A\gamma)^{m+2\alpha m} \bigl( 1 - (B/A)^{m+2\alpha m}
\bigr)\nonumber
\end{eqnarray}
for every even $m\geq2$.
For even $m$ define $G_m$, $H_m$, $k_0$ and $f$ as in the proof of
(b)---Case 2.
If $G_m = \varnothing$ for some even $m$, then $\widetilde{\nu}
(\{ m \}) \leq E_m < 0$ similarly to \textup{(\ref{eq-E2.1})}.
So suppose that $G_m \neq\varnothing$ for all even $m\geq2$.
As seen in the proof of (b), this implies that $G_m$ consists
of a single element and that $f\geq1$. Let
$m=m_j = 2^{jf}$ with $j\in\mathbb{N}$, then the
unique element $(n_j',m_j')$ in $G_{m_j}$
satisfies $m_j'/m_j = c^{jk_0}$. Recall that $m_j'$ is odd by the definition
of $G_m$. For large $j$, we then
have $m_j'/2 > m_j + 2\alpha m_j$, and from \textup{(\ref
{eq-Em-odd})} and
\textup{(\ref{eq-Em-even})} it follows that there exists some constant
$d_2 > 0$ such that
\[
\frac{E_{m_j'}}{|E_{m_j}|} \leq d_2 \sqrt{m_j'}   (A / \gamma
)^{m_j'/2} \to0
\qquad \mbox{as }  j\to\infty,
\]
so that $\nu( \{ m_j \}) \leq E_{m_j} + E_{m_j'} < E_{m_j}/2 <0$ for
large $j$, finishing
the proof of (c).

(d) Suppose $q=0$.
By \textup{(\ref{eq-rho-1b})},
$\rho=\rho_{0,r}$ is Bernoulli distributed with parameter~$r$.
Further, $\mu=\mu_{c,0,r}$ is the distribution
of $\sum_{n=0}^\infty c^{-n} U_n$, where $\{ U_n, n\in\mathbb{N}\}$
is an i.i.d. sequence with
distribution $\rho$. The support of $\mu$ is then a subset of
$[0,c/(c-1)]$. It follows that also $\rho^{\mathrm{sym}}$ and
$\mu^{\mathrm{sym}}$ have bounded support. Moreover none of them
is degenerate to a Dirac measure. Hence they
are not infinitely divisible.
\end{pf*}

\begin{example} \label{ex-1}
(a) Let $p=q>0$. Then $\rho_{q,r}$, $\mu_{c,q,r}$,
$(\rho_{q,r})^{\mathrm{sym}}$ and $(\mu_{c,q,r})^{\mathrm{sym}}$
will all be infinitely divisible if $r\in[0, 3-2\sqrt{2}]$, and
none of them is infinitely divisible if $r>3-2\sqrt{2}
\approx0.17157$. Recall that $r$ is the dependence parameter.

(b)
Let $2p=q>0$. Then $\rho_{q,r}$ and $\mu_{c,q,r}$
will be infinitely
divisible for $r \in[0, (13-3\sqrt{17})/4]$ and fail
to be infinitely divisible for $r > (13-3\sqrt{17})/4 \approx
0.15767$. On the other hand, $(\rho_{q,r})^{\mathrm{sym}}$
and $(\mu_{c,q,r})^{\mathrm{sym}}$ are
infinitely divisible if and only if $r \in[0,(13-3\sqrt{17})/4] \cup
[1/2,1)$.
\end{example}


\section{Continuous-singularity and absolute continuity}\label{sec3}

We continue to study the distribution
\[
\mu_{c,q,r} = \mathcal{L} \biggl( \int_0^{\infty-} c^{-N_{s-}} \, dY_s
\biggr)
\]
defined by a process $\{(N_t,Y_t), t\geq0\}$ and
a constant $c>1$ in Section \ref{sec2}.
The parameters $p$, $q$ and $r$ with $p+q+r=1$ are assumed to satisfy
$p,q,r<1$ and $p,q,r\geq0$ throughout this section (see
the first paragraph of Section \ref{sec2}).
In this section continuity properties of $\mu_{c,q,r}$ are considered.
Since $\mu_{c,q,r}$ is $c^{-1}$-decomposable and nondegenerate,
it is either continuous-singular or absolutely continuous, as
Wolfe's theorem \cite{Wo} says. So our problem is to specify the
continuous-singular case and the absolutely continuous case.
To get complete criteria for the two cases is a difficult problem,
far from being achieved.

We use two classes of numbers,
namely Pisot--Vijayaraghavan (P.V.) numbers (sometimes called Pisot
numbers) and Peres--Solomyak (P.S.) numbers.
A number $c>1$ is called a \textit{P.V. number}
if there exists a polynomial $F(x)$ with integer coefficients
with leading coefficient 1 such that $c$ is a simple root of $F(x)$
and all other roots have
a modulus of less than 1. Every positive integer greater than 1
is a P.V. number, but also $(1+\sqrt{5})/2$ and the unique
real root of $x^3-x-1=0$ are nontrivial
examples. There exist countably infinitely many P.V. numbers
which are not integers.
See Peres, Schlag and Solomyak \cite{PSS} for related information.
On the other hand, following Watanabe \cite{Wa00}, we call $c^{-1}$
a \textit{P.S. number} if $c>1$ and if there are $p_0\in(1/2,1)$ and
$k\in
\mathbb{N}$
such that the $k$th power of the characteristic function of
the distribution of $\sum_{n=0}^{\infty}
c^{-n}U_n$, where $\{U_n\}$ is Bernoulli i.i.d. with $P[U_n=0]=1-P[U_n=1]
=p_0$, is integrable.
Watanabe \cite{Wa00} points out that the paper \cite{PS}
of Peres and Solomyak
contains the proof that the set of P.S. numbers in
the interval $(0,1)$ has
Lebesgue measure $1$. However, according to \cite{Wa00}, an explicit
example of a P.S. number is not known so far. As follows from the
results of \cite{Wa00}, the set of P.V. numbers and the set of reciprocals
of P.S. numbers are disjoint.

\begin{theorem}\label{thm-PV}
Assume that $c$ is a P.V. number and that $q>0$. Then
$\mu_{c,q,r}$ is continuous-singular.
\end{theorem}

Recall that the assumption $q>0$ merely excludes the case of infinite
Bernoulli convolutions.

\begin{pf*}{Proof of Theorem \protect\ref{thm-PV}}
Write $\mu=\mu_{c,q,r}$.
The following proof of continu\-ous-singularity of
$\mu$ is based on an idea of Erd\H{o}s \cite{Er}.
It is enough to show that it is not
absolutely continuous. Thus,
by virtue of the Riemann--Lebesgue theorem,
it is enough to find a sequence $z_k\to\infty$ such
that
\[
\limsup_{k\to\infty} |\widehat\mu(z_k)|>0.
\]
By the definition of a P.V. number, there is a polynomial
$
F(x)=x^N+a_{N-1}x^{N-1}+\cdots+a_1 x+a_0
$
such that $a_{N-1},\ldots,a_0\in\mathbb{Z}$, $F(c)=0$, and the totality
$\{\alpha_1,\ldots,\alpha_N\}$ of roots of $F(x)$ satisfies $\alpha
_1=c$ and
$|\alpha_j|<1$ for $2\leq j\leq N$.
Choose
$z_k=2\pi c^k$.
Now we divide the proof into three cases: (Case 1) $p>0$ and
$r\leq pq$; (Case 2) $p=0$; (Case 3) $p>0$ and
$r> pq$. Recall that $q>0$ is always assumed.

\textit{Case} 1:
As in the proofs of Theorem \ref{thm-id}(b)---Case 2, we have
%
\begin{equation}\label{new1}
\widehat\mu(z)=
\exp\Biggl[\sum_{n=0}^{\infty}
\sum_{m=1}^{\infty}(e^{imc^{-n}z}-1)a_m\Biggr],
\qquad  z\in\mathbb{R},
\end{equation}
with
%
\begin{equation}\label{new2}
a_m = m^{-1} q^m \bigl(1 - \bigl(-r/(pq)\bigr)^m\bigr)\geq0.
\end{equation}
The double series above is absolutely convergent. We then have
%
\begin{equation}\label{new3}
|\widehat\mu(z)|=
\exp\Biggl[-\sum_{n=0}^{\infty}
\sum_{m=1}^{\infty}\bigl(1-\cos(mc^{-n}z)\bigr)a_m\Biggr].
\end{equation}
Thus
\[
|\widehat\mu(z_k)|=\exp\Biggl[-\sum_{n=0}^{\infty}
\sum_{m=1}^{\infty}\bigl(1- \cos(2\pi mc^{k-n})\bigr) a_m
\Biggr]=\exp\Biggl[-\sum_{m=1}^{\infty}(S_m+R_m)
a_m\Biggr]
\]
with
\[
S_m=\sum_{n=0}^k
\bigl(1- \cos(2\pi mc^{k-n})\bigr),\qquad  R_m=\sum_{n=k+1}^\infty
\bigl(1- \cos(2\pi mc^{k-n})\bigr).
\]
Now
\[
S_m=\sum_{n=0}^{k}
\bigl(1-\cos(2\pi mc^{n})\bigr)
=\sum_{n=0}^{k}
\Biggl(1-\cos\Biggl(2\pi m\sum_{j=2}^N {\alpha_j}^n\Biggr)\Biggr),
\]
since
$c^n=\sum_{j=1}^N {\alpha_j}^n-\sum_{j=2}^N {\alpha_j}^n$
and $\sum_{j=1}^N {\alpha_j}^n$ is an integer.
The latter is a consequence of the symmetric function theorem
in algebra (e.g., Lang \cite{La}, Section~IV.6),
implying that $\sum_{j=1}^N {\alpha_j}^n$,
as a symmetric function
of $\alpha_1, \ldots, \alpha_N$,
can be expressed as a polynomial with integer coefficients
in the elementary symmetric functions,
which are integer valued themselves since $F$ has integer
coefficients with leading coefficient 1.
Choose $0 < \delta< 1$ such that $|\alpha_j|<\delta$
for $j=2, \ldots, N$. Then, with some constants
$C_1, C_2,C_3$,
\begin{eqnarray*}
S_m&\leq& C_1\sum_{n=0}^{k}\Biggl(m
\sum_{j=2}^N {\alpha_j}^n\Biggr)^2
\leq C_2 m^2 \sum_{n=0}^{k} \sum_{j=2}^N |{\alpha_j}|^{2n}\\
&\leq& C_3 m^2\sum_{n=0}^{k}\delta^{2n}\leq C_3 m^2/(1-\delta^2).
\end{eqnarray*}
Further, we have
\[
R_m\leq C_1
\sum_{n=1}^\infty(m c^{-n})^2 = C_1 m^2 / (c^2 -1).
\]
Hence, it follows that
\[
|\widehat{\mu}(z_k)| \geq\exp\Biggl[ -
\sum_{m=1}^{\infty} a_m m^2 \biggl(
\frac{C_3}{1-\delta^2} + \frac{C_1}{c^2-1} \biggr) \Biggr].
\]
This shows that
$\limsup_{k\to\infty} |\widehat{\mu}(z_k)| > 0$,
since $\sum_{m=1}^{\infty} a_m m^2<\infty$.

\textit{Case} 2: We have
\[
|\widehat\rho_{q,r}(z)|=\exp\Biggl[ \sum_{m=1}^{\infty} (\cos mz-1)
\frac{q^m}{m}\Biggr]
\]
by the remark at the beginning of the proof of
Lemma \ref{lem-rho-id}(a). Hence the situation is the same as
in Case 1.

\textit{Case} 3: Recall the proof of Lemma \ref{lem-rho-id}(c).
We have
\[
|\widehat\rho_{q,r}(z)|^2=\exp\Biggl[ 2\sum_{m=1}^{\infty}E_m
(\cos mz-1)
\Biggr]
\]
with $E_m$ of \textup{(\ref{eq-Em})}. Hence
\[
|\widehat\mu(z)|^2=\prod_{n=0}^{\infty} |\widehat\rho_{q,r}(c^{-n}z)|^2
\geq\exp\Biggl[ -2\sum_{n=0}^{\infty}\sum_{m=1}^{\infty}
E_m^+ \bigl(1-\cos(mc^{-n}z)\bigr)\Biggr],
\]
where $E_m^+=\max\{E_m,0\}$. We have $\sum_{m=1}^{\infty}
E_m^+ m^2<\infty$, since
\begin{eqnarray*}
\sum_{m=1}^{\infty}
E_m^+ m^2&\leq&\sum_{m=1}^{\infty}\sum_{h=0}^{\infty}
m^2|D_{m+2h}|\pmatrix{m+2h\cr h}\\
&=&\sum_{m=1}^{\infty}\sum_{k-m\ \mathrm{even}\ \geq0} m^2
|D_k|\pmatrix{k\cr(k-m)/2}\\
&\leq&\sum_{k=1}^{\infty} k^2 |D_k|
\sum_{l=0}^{\lfloor(k-1)/2\rfloor} \pmatrix{k\cr l}
\leq\sum_{k=1}^{\infty} \frac{k}{2} (A^k+ B^k)<\infty,
\end{eqnarray*}
noting that $\sum_{l=0}^{\lfloor(k-1)/2\rfloor} {k\choose l}
\leq2^{k-1}$ and $|D_k|\leq k^{-1} 2^{-k} (A^k+B^k)$ with
$0<A<1$ and $0<B<1$. Hence we obtain $\limsup_{k\to\infty}
|\widehat\mu(z_k)|^2>0$ exactly in the same way as in Case 1.\vadjust{\goodbreak}
\end{pf*}

\begin{theorem}\label{thm-PS}
Assume that $c^{-1}$ is a P.S. number. Then there exists
$\varepsilon=\varepsilon(c)\in(0,1)$ such that $\mu_{c,q,r}$ is absolutely
continuous with bounded continuous density whenever $p>0$, $r\leq
pq$ and $q\geq1-\varepsilon$, or whenever $p=0$ and $q\geq
1-\varepsilon$. In
particular, there exist constants $C_1 = C_1(c)>0$ and $C_2 =
C_2(c)>0$ such that $\mu_{c,q,r}$ is absolutely continuous with
bounded continuous density whenever $q \geq C_1 p \geq C_2 r$.
\end{theorem}

Recall that  Lebesgue almost all $c\in(1,\infty)$ are the reciprocals
of P.S. numbers.

\begin{pf*}{Proof of Theorem \protect\ref{thm-PS}}
Let $\mu=\mu_{c,q,r}$. Let $p_0\in(1/2,1)$ and $k\in\mathbb{N}$
as in the definition of a P.S. number.
The following proof was suggested by an argument of
Watanabe \cite{Wa00}.
Let $K := k |\log(2p_0-1)|/2$, which is positive.
Then (2.4) of~\cite{Wa00} tells us that
%
\begin{equation} \label{eq-PS-2}
\int_{-\infty}^{\infty}\exp\Biggl\{ \alpha\sum_{n=0}^{\infty}
\bigl(\cos(c^{-n}u)-1\bigr)\Biggr\} \,du
<\infty \qquad\mbox{whenever }\alpha\geq K.
\end{equation}
Under the condition that $p>0$, $q>0$ and $r\leq pq$, we have
\textup{(\ref{new3})} with $a_m$ of \textup{(\ref{new2})}. Let
$\alpha_0=\sum_{m=1}^{\infty} a_m$. Then it follows from Jensen's
inequality that
\begin{eqnarray*}
\int_{-\infty}^{\infty}|\widehat\mu(z)|\,dz&=&2\int_{0}^{\infty
}\exp
\Biggl[
\frac{1}{\alpha_0} \sum_{m=1}^{\infty} a_m\Biggl(\alpha_0\sum
_{n=0}^{\infty}
\bigl(\cos(mc^{-n}z)-1\bigr)\Biggr)\Biggr] \,dz\\
&\leq&2\int_{0}^{\infty}\Biggl[\frac{1}{\alpha_0} \sum
_{m=1}^{\infty
} a_m
\exp\Biggl( \alpha_0\sum_{n=0}^{\infty}
\bigl(\cos(mc^{-n}z)-1\bigr)\Biggr)\Biggr] \,dz\\
&=&\frac{2}{\alpha_0}\Biggl(\sum_{m=1}^{\infty}\frac{a_m}{m}\Biggr)
\int_0^{\infty}
\exp\Biggl( \alpha_0\sum_{n=0}^{\infty}
\bigl(\cos(c^{-n}u)-1\bigr)\Biggr) \,du.
\end{eqnarray*}
The last integral is finite whenever $\alpha_0 \geq K$ by
\textup{(\ref{eq-PS-2})}. We have $a_m\geq m^{-1}q^m$ for~$m$ odd,
and it
follows that $\alpha_0$ tends to $\infty$ as $q\uparrow1$. Thus
there is
$\varepsilon=\varepsilon(c)$ such that $\alpha_0\geq K$ for all
$q\geq1-\varepsilon$. Hence
$\mu$ has bounded continuous density whenever $p>0$, $r\leq pq$ and
$q\geq1-\varepsilon$. The case when $p=0$ and $q\geq1-\varepsilon$ follows
similarly, with $a_m=m^{-1}q^m$ in the above calculations.

To see the second half of the theorem, suppose that
$q \geq C_1 p \geq C_2 r$ with $C_1,C_2>0$.
Then $p>0$ since $p+r>0$, and hence $q>0$. Thus
\[
q=\biggl(1+\frac{p}{q}+\frac{r}{q}\biggr)^{-1}
\geq\biggl(1+\frac{p}{q}+\frac{C_1 p}{C_2 q}\biggr)^{-1}\geq
\biggl(1+\frac{1}{C_1}+\frac{1}{C_2}\biggr)^{-1}.
\]
Hence, $q\geq1-\varepsilon(c)$ if $C_1$ and $C_2$ are large enough.
We also have
\[
\frac{r}{pq}\leq\frac{C_1}{C_2}\biggl(1+\frac{1}{C_1}+\frac
{1}{C_2}\biggr).
\]
Hence $r/(pq)\leq1$ if $C_1$ is fixed and
$C_2$ is large.
Thus there are $C_1$ and $C_2$ such that $q\geq1-\varepsilon(c)$,
$p>0$ and
$r/(pq)\leq1$ whenever $q\geq C_1 p\geq C_2 r$.
\end{pf*}

Now we use the \textit{entropy} $H(\rho)$ of a discrete
probability measure $\rho$ on $\mathbb{R}$. Here discrete means that
$\rho$ is
concentrated on a countable set. We define
\[
H(\rho) := - \sum_{a\in C} \rho( \{ a \} )
  \log\rho( \{ a \} ),
\]
where $C$ is the carrier (the set of
points with positive mass) of $\rho$.

\begin{theorem}\label{thm-dim} Assume that $q>0$. Then
%
\begin{equation}\label{huvw}\qquad
H(\rho_{q,r})= (q+r) \biggl( \log\frac{1}{1-q}
+ \frac{1}{1-q} \log\frac{1}{q} - \log\frac{q+r}{q}
\biggr) + p \log\frac{1}{p},
\end{equation}
where $p\log(1/p)$ is understood to be zero for $p=0$.
The following are true:

\textup{(a)} The Hausdorff dimension of $\mu_{c,q,r}$ is estimated as
%
\begin{equation}\label{Hdimuvw}
\dim (\mu_{c,q,r})\leq\frac{H(\rho_{q,r})}{\log c}.
\end{equation}

\textup{(b)} For each $c>1$, there exists a constant $C_3 = C_3(c)>0$
such that $\mu_{c,q,r}$ is continuous-singular whenever
$\max\{q,r\} \leq C_3 p$.

\textup{(c)} Fix $q$ and $r$. Then there exists a constant $C_4 = C_4(q,r)>0$
such that $\mu_{c,q,r}$
is continuous-singular whenever $c \geq C_4$.
\end{theorem}

The estimate \textup{(\ref{Hdimuvw})} is meaningful only when
$H(\rho_{q,r})/\log c<1$, as the Hausdorff dimension of any measure
on the line is less than or equal to $1$. In this case~\textup{(\ref{Hdimuvw})}
not only tells the continuous-singularity of
$\mu_{c,q,r}$, but also gives finer information on a set of
full measure for $\mu_{c,q,r}$.

\begin{pf*}{Proof of Theorem \protect\ref{thm-dim}}
Recall that $\rho_{q,r}$ is defined by \textup{(\ref{eq-rho-1})}.
The geometric distribution $\sigma_q$ has entropy
%
\begin{equation}\label{ent-geo}
H(\sigma_q)=-\log(1-q)-\frac{q}{1-q}\log q
\end{equation}
and $H(\rho_{q,r})$ is readily calculated as
\[
H(\rho_{q,r})= ( 1 +r/q ) [ H(\sigma_q) + (1-q) \log(1-q) - q
\log(1+r/q) ] - p \log p,
\]
which shows \textup{(\ref{huvw})}.

\textup{(a)} Applying the remarkable Theorem 2.2 of Watanabe \cite{Wa00}
on $c^{-1}$-decompo\-sable distributions, we obtain \textup{(\ref
{Hdimuvw})}.

\textup{(b)} Notice that, since $\mu$ is continuous-singular or absolutely
continuous, it must be continuous-singular if its
Hausdorff dimension is less than 1.
Suppose that $\max\{q,r\} \leq C_3 p$ with $C_3>0$. Then $p>0$
and hence $p=(1+q/p+r/p)^{-1}\geq(1+2C_3)^{-1}$,
which tends to 1 as $C_3\to0$.
Hence $q\to0$ and $r\to0$ as $C_3\to0$.
Thus
%
\[
H(\rho_{q,r})=(q+r)\biggl( \log\frac{1}{1-q} + \frac{q}{1-q} \log
\frac{1}{q} + \log\frac{1}{q+r}\biggr) +p \log\frac{1}{p} \to0
\]
and hence
\[
\sup\bigl\{ H(\rho_{q,r})\dvtx\max\{q,r\}\leq C_3 p\bigr\}\to0,\qquad
C_3\to0.
\]
This shows (b).

\textup{(c)} For given $q,r$, take any $C_4>\exp H(\rho_{q,r})$. Then the
assertion follows from the estimate \textup{(\ref{Hdimuvw})}.
\end{pf*}

It is known that the distributions of some L\'evy processes have
time evolution, that is, change their qualitative properties as
time passes (see Chapters 5 and 10 of Sato \cite{Sa} and Watanabe
\cite{Wa01}). It was Watanabe \cite{Wa00} who showed that
the distributions of some
semi-self-decomposable processes have time evolution in
continuous-singularity and absolute continuity.
It is of interest that the L\'evy process $\{Z_t,t\geq0\}$ determined
by the distribution $\mu_{c,q,r}$ when it is infinitely divisible
gives an explicit example of time evolution of this sort as the following
theorem shows. Note that $\mathcal L(Z_t)=(\mu_{c,q,r})^{t*}$ for
$t\geq0$.

\begin{theorem} \label{thm-power-2}
Assume that either $p>0$, $q>0$ and $r\leq pq$ or $p=0$.
Write $\rho=\rho_{q,r}$ and $\mu=\mu_{c,q,r}$.
Then the following are true:

\textup{(a)} There are $t_1=t_1(c,q,r)$ and
$t_2=t_2(c,q,r)$ with $0<t_1\leq t_2\leq\infty$ such that ${\mu
}^{t*}$ is
continuous-singular for all $t\in(0,t_1)$,
absolutely continuous without bounded continuous
density for all $t\in(t_1, t_2)$ if
$t_1 < t_2$, and
absolutely continuous with bounded continuous
density for all $t\in(t_2,\infty)$ if $t_2<\infty$.

\textup{(b)} If $c$ is a P.V. number, then $t_1=\infty$, that is,
${\mu}^{t*}$ is continuous-singular for all $t>0$.

\textup{(c)} If $c^{-1}$ is a P.S. number, then $t_2<\infty$. Thus
$t_2<\infty$
for Lebesgue almost all $c>1$.

\textup{(d)} The entropy $H(\rho^{t*})$ is a finite, continuous, strictly
increasing function
of $t\in[0,\infty)$, vanishing at $t=0$, and
%
\begin{equation}\label{Hdim}
\dim ({\mu}^{t*})\leq\frac{H(\rho^{t\ast})}{\log c} \qquad\mbox
{for all
$t\geq0$}.
\end{equation}

\textup{(e)} Fix $q$, $r$ and $t>0$. Then
${\mu}^{t*}$ is continuous-singular for all
sufficiently large $c>1$.
\end{theorem}

\begin{lemma} \label{lem-ent}
If $\rho$ is a distribution on $\mathbb{Z}$ with finite absolute
moment of order $1+\varepsilon$ for some $\varepsilon>0$,
then its entropy $H(\rho)$ is finite.
\end{lemma}

\begin{pf}
Let $\rho=\sum_{m=-\infty}^{\infty} p_m \delta_m$.
Then $\sum_{m=-\infty}^{\infty} |m|^{1+\varepsilon} p_m<\infty$. Hence
there is a
constant $C>0$ such that $p_m\leq C|m|^{-1-\varepsilon}$. The function
$f(x)=x\log(1/x)$ is increasing for $0\leq x\leq e^{-1}$.
Hence,
\[
H(\rho)\leq\sum_{|m|\leq m_0}p_m \log(1/p_m) +
\sum_{|m|> m_0} C|m|^{-1-\varepsilon} \log((C|m|^{-1-\varepsilon
})^{-1})<\infty
\]
with an appropriate choice of $m_0$.
\end{pf}

\begin{pf*}{Proof of Theorem \protect\ref{thm-power-2}}
Observe that under our assumption $\mu$ is infinitely divisible
(see Theorem \ref{thm-id}),
so that $\mu^{t*}$ is definable for all $t\geq0$.
If $p>0$ and $r\leq pq$, it follows from \textup{(\ref{eq-cf2})}
that $\rho$ is a compound Poisson distribution, concentrated on
$\mathbb{N}_0$,
with finite second moment, since $\sum_{m=1}^{\infty}m^2 a_m<\infty$
for $a_m = \nu_\rho(\{ m \})$. If $p=0$, then $\rho$ is a geometric
distribution shifted by 1 (see the proof of Lemma \ref{lem-rho-id}).
In both cases, $H(\rho) < \infty$ by
Lemma \ref{lem-ent}.
The property that $\mu$ is $c^{-1}$-decomposable
is preserved to convolution powers, since \textup{(\ref{eq-SS1})} implies
\[
\widehat{{\mu}^{t*}}(z) = \widehat{\rho^{t*}} (z)   \widehat{{\mu}^{t*}}
(c^{-1} z)
\]
for any $t\ge0$. Thus we have
%
\begin{equation} \label{extra-5}\quad
\widehat{\mu^{t*}} (z) = \prod_{n=0}^\infty\widehat{\rho^{t*}}
(c^{-n} z) =
\exp\Biggl( it \gamma_p^0 z + t
\sum_{n=0}^\infty\sum_{m=1}^\infty(e^{im c^{-n} z} -1) a_m
\Biggr),
\end{equation}
where $\gamma_p^0=0$ for $p>0$
and $\gamma_p^0 = \sum_{n=0}^\infty c^{-n} = c/(c-1)$ for $p=0$.
If $p>0$, $q>0$ and $r=0$, then $\rho$ is geometric.
If $p=0$, then the result is
the same as in the case of $\rho$ being a geometric distribution, since
shifts do not change entropy, Hausdorff dimension and continuity
properties.
Hence, from now on, we exclude the case $p=0$.
Let us prove (d) first.
We obtain the properties of $H(\rho^t)$ from Proposition~5.1 of
Watanabe \cite{Wa00} or Exercise 29.24 of Sato \cite{Sa}.
Applying
Theorem 2.2 of \cite{Wa00} to the law ${\mu}^{t*}$, we get \textup
{(\ref{Hdim})}.

(a) It follows from (d) that ${\mu}^{t*}$ is
continuous-singular for all sufficiently
small $t>0$. Let $t_1$ be the supremum of $t>0$ for which $\mu^{t*}$ is
continuous-singular. Then $0<t_1\leq\infty$. Noting that if $t<t'$, then
$\mu^{t*}$ is a convolution factor of $\mu^{t'*}$, use Lemma 27.1 of
\cite{Sa}. Thus ${\mu}^{t*}$ is continuous-singular for all
$t\in(0,t_1)$ and absolutely continuous for all
$t\in(t_1,\infty)$ if $t_1<\infty$.
Let $t_2$ be the infimum of $t>t_1$ for which $\mu^{t*}$
has bounded continuous density.
Use the fact that a distribution is absolutely
continuous with bounded continuous density if a convolution factor of
it has this property. Then we obtain the assertion (a).

(b) The proof is the same as that of Theorem
\ref{thm-PV} with $a_m$ replaced by $ta_m$.

(c) Let $c^{-1}$ be a P.S. number. Then ${\mu}^{t*}$ is absolutely
continuous with bounded continuous density for all sufficiently large
$t$, which is shown in the same way as
Theorem \ref{thm-PS}, or we can apply
Theorem 2.1 of Watanabe \cite{Wa00}.

(e) Obvious consequence of (d).
\end{pf*}

In the rest of this section we consider the case where $\{N_t\}$ and
$\{Y_t\}$ are independent, that is, the case where
$r=0$. The following theorem is
largely a repetition of Theorems \ref{thm-PV}--\ref{thm-dim} in this case
but, since the L\'evy measure of $\mu_{c,q,0}$ is increasing with
respect to $q$, we obtain stronger statements.

\begin{theorem} \label{thm-t-fixed}
Under the condition that $0<q<1$ and $r=0$ (i.e., $\{N_t\}$ and~$\{Y_t\}$ are independent), the following are true:

\textup{(a)} For any $c>1$ there are constants $q_1=q_1(c)$ and $q_2=q_2(c)$
satisfying
$0 < q_1 \leq q_2 \leq1$ with the following properties: $\mu_{c,q,0}$ is
continuous-singular for all $q\in
(0, q_1)$, absolutely continuous without bounded continuous
density for all $q \in(q_1, q_2)$ if
$q_1 < q_2$, and absolutely continuous with bounded
continuous density for all $q \in(q_2, 1)$ if $q_2<1$.

\textup{(b)} If $c$ is a P.V. number, then $q_1=1$, that is, $\mu_{c,q,0}$ is
continuous-singular for all $q$.

\textup{(c)} If $c^{-1}$ is a P.S. number, then $q_2<1$. Hence $q_2<1$ for
Lebesgue almost all $c>1$.

\textup{(d)} The Hausdorff dimension of $\mu_{c,q,0}$ is estimated as
%
\begin{equation}\label{Hdimbeta}
\dim (\mu_{c,q,0})\leq\frac{H(\rho_{q,0})}{\log c},
\end{equation}
where
%
\begin{equation} \label{h-fun}
H(\rho_{q,0})=-\log(1-q)- \frac{q}{1-q} \log q,
\end{equation}
which is a finite, continuous, strictly increasing function of
$q\in(0,1)$ and tends to~$0$ as $q\downarrow0$.

\textup{(e)} Fix $c>1$. If $q$ is so small that $H(\rho_{q,0}) < \log c$,
then $\mu_{c,q,0}$ is continuous-singular. In particular, if
%
\begin{equation}\label{smallq}
0 < q< 1 -(\log2/\log c),
\end{equation}
then $\mu_{c,q,0}$ is continuous-singular.

\textup{(f)} Fix $q >0$. If $c$ is so large that $c>\exp H(\rho_{q,0})$,
then $\mu_{c,q,0}$ is continuous-singular.
\end{theorem}

\begin{pf}
Let us begin with the proof of (d). The estimate \textup{(\ref{Hdimbeta})}
follows from~\textup{(\ref{Hdimuvw})} of Theorem \ref{thm-dim}. The
expression
\textup{(\ref{h-fun})} is exactly \textup{(\ref{ent-geo})}, since
$\rho_{q,0}=\sigma_q$.

(e) and (f) These come from (d), as a distribution with
Hausdorff dimension $<1$ cannot be absolutely continuous. We get
the sufficient condition \textup{(\ref{smallq})}, since
\[
H(\rho_{q,0}) = \frac{1}{1-q} \biggl( (1-q)\log\frac{1}{1-q}+q
\log\frac{1}{q} \biggr) \leq \frac{1}{1-q}\log2
\]
by strict concavity of the function $\log x$.

(a) Recall that $\mu_{c,q,0}$ has L\'evy measure
\[
\nu_{c,q,0}=\sum_{n=0}^{\infty}
\sum_{m=1}^{\infty} \frac{q^m}{m}
\delta_{c^{-n}m}.
\]
Hence, if $q<q'$, then $\mu_{c,q,0}$ is a convolution factor
of $\mu_{c,q',0}$. Now the proof is obtained by the same argument as
in the proof of (a) of
Theorem \ref{thm-power-2}.

(b) Consequence of Theorem \ref{thm-PV}.

(c) This follows from Theorem \ref{thm-PS}.
\end{pf}

\begin{example}\label{newex}
(a) In the case $c=e$,
$\mu_{e,q,0}$ is continuous-singular if
\[
q\leq1-\log2\approx0.30685.
\]
This follows from \textup{(\ref{smallq})} in Theorem \ref
{thm-t-fixed}.

(b)
In the case $q=1/2$, $\mu_{c,1/2,0}$
is continuous-singular if $c>4$ since $H(\rho_{1/2,0})=2\log2$,
as (f) of Theorem \ref{thm-t-fixed}
says.
\end{example}

In the independent case ($0<q<1$ and $r=0$), the assumption in
Theorem \ref{thm-power-2} is satisfied. So the assertions
on time evolution of $\mu_{c,q,0}$ hold true as in that theorem.
It is of interest to estimate $H({\rho_{q,0}}^{t*})$ appearing
in the right-hand side of \textup{(\ref{Hdim})}.

\begin{proposition} \label{prop-2C}
If $0<q<1$, then
%
\begin{equation} \label{g-fun}
H({\rho_{q,0}}^{t*})\leq t \biggl[\frac{1}{p}\biggl( 1 + 2 \log\frac{1}{p}
\biggr)
+\frac{q}{p} \log\frac{1}{t}\biggr] \qquad\mbox{for } 0 < t
\leq1,
\end{equation}
where $p=1-q$. The right-hand side of \textup{(\ref{g-fun})} is a strictly
increasing function of $t\in(0,1]$ which tends to $0$ as $t\downarrow0$.
\end{proposition}

\begin{pf}
Write $\rho=\rho_{q,0}$. Since $\rho$ equals geometric distribution
$\sigma_q$ with parameter $p$, the
distribution $\rho^{t*}$, $t>0$, is a negative binomial
distribution with
parameters $t$ and $p$, that is,
\[
\rho^{t*} ( \{ k \} )= \pmatrix{ -t \cr k } p^t (-q)^k,
 \qquad k\in\mathbb{N}_0.
\]
To estimate $H(\rho^{t*})$ from above,
observe that $t p^t q^k/k \leq\rho^{t*} (\{ k \} )
\leq t   q^k$ for $0 < t \leq1$ and $k \in\mathbb{N}$
so that
\begin{eqnarray*}
H(\rho^{t*}) & = &- \sum_{k=0}^\infty\rho^{t*}
(\{ k \})   \log\rho^{t*} ( \{ k \})\\
&\leq&- (\log p^t) + \sum_{k=1}^\infty t q^k
\bigl( \log k - \log(t p^t) - k \log q \bigr) \\
& \leq& t \biggl[ \log\frac{1}{p} + \frac{1}{p} \log\frac{1}{p}
- \frac{q}{p} \log(t p^t) - \frac{q}{p^2} \log q \biggr],
\end{eqnarray*}
where we used $\sum_{k=1}^\infty k   q^k = q/p^2$ and
\[
\sum_{k=1}^\infty q^k \log k \leq\sum_{k=1}^\infty q^k
\sum_{n=1}^k \frac{1}{n} =
\frac{1}{p} \log\frac{1}{p},
\]
compare Gradshteyn and Ryzhik \cite{GR}, Formula 1.513.6.
Recalling that $p^t \geq p$ since $t\leq1$ and this can be further
estimated to
\[
H(\rho^{t*}) \leq t \biggl[ \frac{2}{p}\log\frac{1}{p}
+ \frac{q}{p}\log\frac{1}{t} +
\frac{q}{p^2} \log\frac{1}{q} \biggr].
\]
Together with $(q/p) \log(1/q)=(q/p)\log(1 + p/q) \leq1$,
this gives \textup{(\ref{g-fun})}.
\end{pf}

\section*{Acknowledgments}
We thank Toshiro Watanabe for
drawing our attention to~\cite{PSS} and giving information on
infinite Bernoulli convolutions.


%
\printaddresses

\end{document}